\documentclass{amsart}
\theoremstyle{plain}

\errorcontextlines=0 \numberwithin{equation}{section}

\begin{document}
\large
\title[B-minimal sub-manifolds and their stability]
{B-sub-manifolds and their stability}
\author{Li MA }

\address{Department of mathematical sciences \\
Tsinghua university \\
Beijing 100084 \\
China}

\email{lma@math.tsinghua.edu.cn} \dedicatory{}
\date{Oct. 15th, 2002}


\begin{abstract}
In this paper, we introduce a concept of B-minimal sub-manifolds
and discuss the stability of such a sub-manifold in a Riemannian
manifold $(M,g)$. Assume $B(x)$ is a smooth function on $M$. By
definition, we call a sub-manifold $\Sigma$ {\em B-minimal} in
$(M,g)$ if the product sub-manifold $\Sigma\times S^1$ is a {\em
minimal} sub-manifold in a warped product Riemannian manifold
$(M\times S^1, g+e^{2B(x)}dt^2)$,  so its stability is closely
related to the stability of solitons of mean curvature flows as
noted earlier by G. Huisken , S. Angenent, and K. Smoczyk. We can
show that the "grim reaper" in the curve-shortening problem is
stable in the sense of "symmetric stable" defined by K. Smoczyk.
We also discuss the graphic B-minimal sub-manifold in $R^{n+k}$.
\end{abstract}
\maketitle

\section{ Introduction}

  The aim of this paper is two-folds. One is to introduce a new concept of a
  stationary sub-manifold, which will be called the $B$-minimal
  sub-manifold. The other is to answer a question posed by K.
  Smoczyk [10].

Let $(M^m,g)$ be a Riemannian manifold. Assume that $B$ is a
smooth function on $M^m$. Let $N=M\times S^1$ and let
$ds^2=g+e^{2B}dt^2$, where $t\in S^1$.

{\bf Definition 1:} {\em Let $\Sigma^k$ be a $k$-dimensional
manifold and Let $F:\Sigma\to M$ be a immersion. We say that a
$k$-dimensional sub-manifold $\Sigma^k$ is a $B$-minimal
sub-manifold in $M$ if the immersion
$$(x,t)\in \Sigma\times S^1\to (F(x),t)\in N$$ is a
minimal sub-manifold in $(N, ds^2=g+e^{2B}dt^2)$.}

Assume that a $k$-dimensional sub-manifold $\Sigma^k$ is a
$B$-minimal sub-manifold in $M$. Note that the volume of
$\Sigma\times S^1$ in $N$ is
$$
V(\Sigma)=2\pi\int_{\Sigma}e^B(x)dv_g
$$
where $dv_g$ is the induced volume form in $\Sigma^k$. Hence, a
$B$-minimal sub-manifold is in fact a minimal sub-manifold in $M$
equipped with the conformal metric $e^{2B/k}g$.

{\bf Definition 2:}  {\em We say that a $B$-minimal sub-manifold
in $M$ is stable if if the second variational derivative of the
volume functional $V(\Sigma)$ is positive semi-definite.}

{\bf Remark 3:} {\em Let $X=grad B$ be the vector field on $M$. If
$X$ is a Killing or a conformal vector field, Then $B$-minimal
submanifolds are the soliton solutions of the Mean curvature Flow
in $M$ as noted earlier by G. Huisken [7], S. Angenent [1], and K.
Smoczyk [10]. With this understanding, we can easily get the first
and second variational formulae for a $B-$minimal sub-manifold
(see [8][9]).}

Our main results are the following

 {\bf Theorem 1}: {\em Assume that a $k$-dimensional sub-manifold
$\Sigma^k$ is a $B$-minimal sub-manifold in $M$. Then it satisfies
the B-minimal system:
$$
H=(DB)^N.
$$
where $H$ is the mean curvature vector of $\Sigma$ in $M$ and
$(DB)^N$ is the normal part of the derivative $DB$ on $\Sigma$.}

{\bf Theorem 2}: {\em Assume $B(x,y)=y$ in the $xy$-plane. Then
the grim reaper soliton $y=-\log\cos x$ in the curve-shortening
flow is a $B$-minimal sub-manifold of dimension one, which is
stable in both our sense above and the symmetric stability defined
by K. Smoczyk [10]. In fact, we have the following inequality
$$ \int_{\pi/2}^{\pi/2}
\frac{3\cos^2x-1}{4\cos x} u^2(x)dx\leq (u'(x))^2 \cos xdx
$$
where $u\in C^{\infty}_0(-\frac{\pi}{2}, \frac{\pi}{2})$.}

We remark that the inequality above was posed by  K. Smoczyk in
[10].

Theorem 1 is proved in section two. Theorem 2 is proved in section
three. In the last section , We will discuss the graphic B-minimal
sub-manifold in $R^{n+k}$ and some related questions.

\section{B-minimal sub-manifolds}

In this section, we review some useful formulae for minimal
sub-manifolds (see Chern etc.[2] and Yau [11]) and derive the
$B$-minimal sub-manifold system.

Let $(N,ds^2)$ be a Riemannian manifold of dimension $n+d$. We
will use the following ranges of indices.

\begin{align*}
1&\leq A,B, C,D,\cdots \leq n+d \\
1&\leq \alpha, \beta ,\cdots\leq d\\
1&\leq i,j,k,l,\cdots \leq n.
\end{align*}

Let $(\theta^A)$ be an orthogonal frame on $N$. Let $(e_A)$ be its
dual frame. Write
$$
ds^2=g_{AB}\theta^A\theta^B
$$
and let
$$
(g^{AB})=(g_{AB})^{-1}.
$$

 Then we have the following structure equation:

\begin{align*}
d\theta^{A}& =   \theta^{B}\wedge w_{B}^{A}\\
w_{B}^{A}&=-w^{B}_{A}\\
\Phi^A_B&= dw^A_B-w_{B}^{C}\wedge w_{C}^{A}\\
&=\frac{1}{2}K^A_{BCD}\theta^C\wedge\theta^D \\
K^A_{BCD}&=-K^A_{BDC}.
\end{align*}
Let
\begin{align*}
\theta_A&=\theta^A,\\
w_{AB}&=g_{EB}w^E_A, \\
\Phi_{AB}&=g_{EB}\Phi^E_A,\\
K_{BACD}&=g_{EA}K^E_{BCD}.
\end{align*}

Then we have
\begin{align*}
d\theta_{A}& =   \theta_{B}\wedge w_{BA}, \\
w_{BA}&=-w_{AB}, \\
\Phi_{AB}&= dw_{AB}-w_{AC}\wedge w_{CB}, \\
&=\frac{1}{2}K_{ABCD}\theta_C\wedge\theta_D \\
K_{ABCD}&=-K_{BADC}.
\end{align*}

Let $\Sigma$ be an n-dimensional sub-manifold of  $N$. Take an
orthonormal frame $(\theta^A)$ such that $(\theta^i)$ is a local
frame on $\Sigma$ and

$$
\theta^{\alpha}=0, \alpha=1,\cdots,d.
$$
on $\Sigma$. Then we get

$$
d\theta^{\alpha}=0,
$$
that is,
$$
w^{\alpha}_i\wedge\theta^i=0,
$$

By Cartan's lemma we have
$$
w^{\alpha}_i=h^{\alpha}_{ij}\theta^j,
h^{\alpha}_{ij}=h^{\alpha}_{ji}.
$$

The Mean Curvature Vector of $\Sigma$ in $N$ is defined by
$$
H=g^{ij}h^{\alpha}_{ji}e_{\alpha}.
$$

 On the other hand, using the induced
metric on $\Sigma$, we have, on $\Sigma$, the following structure
equations:

\begin{align*}
d\theta^{i}& =    \theta^{j}\wedge w_{j}^{i}\\
w_{j}^{i}&=-w^{j}_{i}\\
\Omega^j_i&= dw^j_i-w_{i}^{k}\wedge w_{k}^{j}\\
&=\frac{1}{2}R^j_{ikl}\theta^k\wedge\theta^l \\
R^j_{ikl}&=-R^j_{ilk}.
\end{align*}

Let

\begin{align*}
\Omega_{ij}&=g_{mj}\Omega^m_i,\\
R_{jikl}&=g_{mi}R^m_{jkl}.
\end{align*}
Then we have , on $\Sigma$,

\begin{align*}
d\theta_{i}& =   \theta_{j}\wedge w_{ji}, \\
w_{ij}&=-w_{ji}, \\
\Omega_{ij}&= dw_{ij}-w_{im}\wedge w_{mj}, \\
&=\frac{1}{2}K_{ijkl}\theta_k\wedge\theta_l \\
K_{ijkl}&=-K_{jikl}.
\end{align*}

Using the structure equations on $N$, we get
$$
dw^i_j-w_{j}^{k}\wedge w_{k}^{i}=w_{j}^{\alpha}\wedge
w_{\alpha}^{i}+\Phi^i_j
$$

Hence, we have
$$
\frac{1}{2}R^i_{jkl}\theta^k\wedge\theta^l=
-h^{\alpha}_{jk}h^{\alpha}_{il}\theta^k \wedge \theta^l+
\frac{1}{2}K^i_{jkl}\theta^k\wedge\theta^l.
$$
and
$$
\frac{1}{2}R_{jikl}\theta^k\wedge\theta^l=
-h^{\alpha}_{jk}h^{\alpha}_{il}\theta^k \wedge \theta^l+
\frac{1}{2}K_{jikl}\theta^k\wedge\theta^l.
$$
Therefore, we obtain the following Gauss equation:
$$
R_{jikl}=K_{jikl}+h^{\alpha}_{ik}h^{\alpha}_{jl}-h^{\alpha}_{il}h^{\alpha}_{jk}
$$

Assume that a $k$-dimensional sub-manifold $\Sigma^k$ is a
$B$-minimal sub-manifold in $M$. Let
$$
N:=M\times S^1
$$
be equipped
with the metric
$$
ds^2=g+e^{2B}dt^2.
$$
Recall that the volume of $\Sigma\times S^1$ in $N$ is
$$
V=2\pi\int_{\Sigma}e^B(x)dv_g
$$
where $dv_g$ is the induced volume form in $\Sigma^k$. Hence, a
$B$-minimal sub-manifold is in fact a minimal sub-manifold in $M$
equipped with the conformal metric $e^{2B/k}g$. We will use this
point of view to get the $B$-minimal sub-manifold equation.

 In the
following, we let $f=B/k$. Let $$\bar \theta^A=e^f\theta^A$$ and
let $$\bar e_A=e^{-f}e_A. $$ Then we have

\begin{align*}
d\bar \theta^A&=e^f(df\wedge \theta^A+d\theta^A)\\
&=\bar\theta^B\wedge (f_B\theta^A-f_A\theta^B+w^A_B)\\
&=\bar\theta^B\wedge \bar w^A_B.
\end{align*}
Here, we used
$$
df=f_A\theta^A.
$$
Hence, we have
$$
\bar w^A_B=f_B\theta^A-f_A\theta^B+w^A_B.
$$

Note that
$$
 w^{\alpha}_j=h^{\alpha}_{jl}
\theta^l
$$
and
\begin{align*}
\bar w^{\alpha}_j&=\bar h^{\alpha}_{jl}\bar \theta^l \\
&=e^f\bar h^{\alpha}_{jl} \theta^l\\
&=-f_{\alpha}\theta^j+w^{\alpha}_j.
\end{align*}
Then we have
$$
-f_{\alpha}\delta_{jl}+h^{\alpha}_{jl}=e^f\bar h^{\alpha}_{jl}.
$$
By this we obtain that
$$
\bar h^{\alpha}_{jl}=e^f(-f_{\alpha}\delta_{jl}+h^{\alpha}_{jl})
$$
and
$$
\bar H=\bar h^{\alpha}_{jj}\bar
e_{alpha}=e^{-2f}(H-kf_{\alpha}e_{\alpha})=e^{-2f}(H-k(Df)^N)
$$
where $(Df)^N$ is the normal part of the derivative $Df$ on
$\Sigma$. Hence the $B$-minimal sub-manifold system is
$$
H=(DB)^N.
$$

\section{A question from K. Smoczyk}

  In his interesting paper [10], K.Smoczyk asks if the following
  inequality is true or not (see (4.21) in [10])(*):
$$
\int_{\pi/2}^{\pi/2} \frac{3\cos^2x-1}{4\cos x} u^2(x)dx\leq
(u'(x))^2 \cos xdx
$$
where $u\in C^{\infty}_0(-\frac{\pi}{2}, \frac{\pi}{2})$. If this
is true, then the grim reaper is symmetric stable defined by K.
Smoczyk. We will prove this inequality. But first of all, let's
see why this inequality is true in another way.

Let $\gamma(x)=(x,y(x))$ be a smooth curve in the plane $R^2$. Let
$B(x)$ be a smooth function in $R^2$.  Then we have
$$
\gamma'(x)=(1,y')
$$
and the unit tangent of the curve is
$$
T=\gamma'/w(x)
$$
where
$$
w(x)=\sqrt{1+{y'}^2}.
$$
The unit normal vector is
$$
N:=(-y',1)/w(x).
$$
The length functional is
$$
L(y)=\int_a^bw(x)dx.
$$
Now we define a new functional
$$
I(y)=\int_a^bw(x)e^{y(x)}dx.
$$
It is easy to see that the first variational formula of $I(\cdot)$
is
\begin{align*}
 \delta I(y)\xi&= \frac{d}{dt}(y+t\xi)_{t=0}\\
 &=\int_{[a,b]}\{\frac{y'\xi'}{w(x)}+w(x)\xi(x)\}e^ydx.
\end{align*}
Here $ \xi\in C^1_0[a,b]. $ Hence a critical point $y=y(x)$ of
$I(\cdot)$ satisfies
$$
-e^{-y(x)}\frac{d}{dx}(e^{y(x)}\frac{y'}{\sqrt{1+{y'}^2}})+\sqrt{1+{y'}^2}=0.
$$

Now we compute directly the second variational formula of
$I(\cdot)$ in the following way.
\begin{align*}
 \delta^ I(y)(\xi, \xi)&= \frac{d}{dt}\delta I(y+t\xi)\xi|_{t=0}\\
&=\int_a^b\{\frac{\xi'}{w(x)}-\frac{{y'}^2{\xi'}^2}{w(x)^3}+2\frac{y'{\xi}'\xi}{w(x)}+
w(x)\xi(x)^2\} e^y\\
&=\int_a^b\{\frac{\xi'}{(w(x)}-\frac{{y'}^2{\xi'}^2}{w(x)^3}\}e^ydx\\
&+\int_a^b\{\frac{y'{(\xi^2)}'}{w(x)}+
w(x)\xi(x)^2\} e^y\\
&=\int_a^b\{\frac{\xi'}{(w(x)}-\frac{y'^2{\xi'}^2}{w(x)^3}\}e^ydx\\
&=\int_a^b\{\frac{{\xi'}^2}{w(x)^3}\}e^ydx\geq 0.
\end{align*}

Hence, any critical point $y=y(x)$ of $I(\cdot)$ is stable in the
sense that the second variational derivative at $y$ is
semi-positive definite.

Let's consider an example. Take
$$
y(x)=-\log \cos x, \quad{x\in (-\pi/2,\pi/2)}.
$$
The curve $(x,y(x))$ is called the grim reaper in the
curve-shortening problem in the plane. Compute
\begin{align*}
y' &=\tan x,\\
e^{y(x)}&=\cos x\\
\sqrt{1+{y'}^2}&=\sqrt{1+{\tan x}^2}=\frac{1}{\cos x}.
\end{align*}
Then one has
\begin{align*}
-e^{-y(x)}\frac{d}{dx}(e^y(x)\frac{y'}{\sqrt{1+{y'}^2}})&= -\cos
x(\cos x \tan x \times \frac{1}{\cos x})'\\
&=-\frac{1}{\cos x}\\
&=-\sqrt{1+{y'}^2}\\
&=\sqrt{1+{\tan x}^2}
\end{align*}
Let $B(x,y)=y$. By definition, this grim reaper is the B-minimal
curve in $R^2$. Therefore, it is stable in our sense above.

Now let's prove the inequality (*) directly. Let $\epsilon>0$ be a
small positive number. We let
$$
p=p_{\epsilon}(x)=\epsilon+\cos x.
$$
and let $J= (-\frac{\pi}{2}, \frac{\pi}{2})$. We define a new
measure
$$
d\mu=p(x)dx
$$
and a new function
$$
f=\frac{3p^2-1}{4p^2}.
$$

 Then we only need to prove the following inequality:
$$
\int_J fu^2d\mu\leq\int_J (u')^2d\mu.
$$
If this is true, we just send $\epsilon \to 0^+$ and get the
inequality (*).

We look for a function $\phi=\phi(x)$ such that
$$\int_J(u'-u\phi)^2d\mu\leq\int_J ((u')^2-fu^2)d\mu.
$$
This is equivalent to
$$
\int_J(u^2(\phi^2+f)-2uu'\phi)d\mu\leq 0.
$$
Note, by using integration by part,
$$
\int_J(2uu'\phi)d\mu=-\int u^2(\phi'-\phi\frac{\sin x}{p})d\mu.
$$
Hence the inequality above can be written as
$$
\int_Ju^2[\phi^2+f+\phi'-\phi\frac{\sin x}{p}]d\mu\leq 0.
$$

We now try to solve the following equation for $\phi$:
$$
\phi'+\phi^2-\phi\frac{\sin x}{p}+f=0
$$
Let $\phi=(\log v)'$. Then we have the following equivalent
equation (*)':
$$
{v''}-\frac{\sin x}{p}v'+f(x)v=0.
$$
Imposing the initial conditions:
$$
v(-\frac{\pi}{2})=1, v'(-\frac{\pi}{2})=0,
$$
we can solve the equation (*)' (see Lemma 1.1 in Chapter IV in P.
Hartman [6]). Hence the inequality (*) is true and we proved
Theorem 2.

\section{Graphic B-Minimal Sub-manifolds}

Let $D\subset R^n$ be a domain of $R^n$. Write $x=(x^1,...,
x^n)\in D$. Let $y=y(x)\in R^k$ be a vector-valued smooth
function. Define the graphic sub-manifold
$$
\Sigma=\{(x,y(x));x\in D\}\subset R^{n+k}.
$$

Let $F(x)=(x,y(x))$ be the graph-mapping in the space $R^n\times
R^k$. Let $B(y)$ be a smooth function in $R^k$.  Then we have
$$
DF(x)=(id,D_xy)
$$
and Let
$$
w(x)=\sqrt{\det(\delta_{ij}+D_{x^i}yD_{x^j}y)}.
$$

Write
$$
g_{ij}=\delta_{ij}+D_{x^i}yD_{x^j}y
$$
and
$$
(g^{ij})=(g_{ij})^{-1}.
$$

 Now we define a new functional
$$
I(y)=\int_Dw(x)e^{B(y(x))}dx.
$$
It is easy to see that the first variational formula of $I(\cdot)$
is
\begin{align*}
 \delta I(y)\xi&= \frac{d}{dt}(y+t\xi)_{t=0}\\
 &=\int_{D}\{\frac{g^{ij}\langle D_{x^i}y,D_{x^j}\xi\rangle}{w(x)}
 +w(x)\langle D_yB,\xi(x)\rangle\}e^{B(y)}dx.
\end{align*}
Here
$
\xi\in C^1_0(D).
$
Hence, $\Sigma$ is a B-minimal sub-manifold in $R^{n+k}$ if and
only if $y=y(x)$ is a critical point of $I(\cdot)$. So $y=y(x)$
satisfies
$$
-e^{-B(y)}D_{x^j}(e^{B(y)}\frac{g^{ij}D_{x^i} y
}{w(x)})+{w(x)}D_yB(y)=0.
$$
In  the special case where $k=1$ and $B(y)=y$, we have
$$
w(x)=\sqrt{1+{|\bigtriangledown y |}^2},
$$
and the B-minimal equation is:
$$
-e^{-y}div(e^y\frac{\bigtriangledown y
}{\sqrt{1+{|\bigtriangledown y|}^2}})+\sqrt{1+{|\bigtriangledown y
|}^2}=0.
$$

Then one may study the Bernstein problem or a-priori estimates for
the graphic B-minimal sub-manifolds. One can also propose a
similar concept like $B$-harmonic maps between Riemannian
manifolds $(M,g)$ and $(N,h)$ with a given smooth function $B$
defined on $N$ in the following way: Let $u:M\to N$ be a $C^1$
mapping and let
$$
e(u)=|du|^2e^{B(u)}
$$
be its energy density on $M$. Define the B-energy functional as
$$
E(u)=\frac{1}{2}\int_M e(u)dx
$$
Then we call $u$ a {\em B-harmonic map} if it is a critical point
of the $B$-energy functional. One may discuss the Liouville
property for B-Harmonic maps. But we will not discuss this concept
more in this paper .

\end{document}